\documentclass[a4paper,11pt]{amsart}
\usepackage{amsmath}
\usepackage{amsaddr}
\usepackage{amsthm}
\usepackage{graphicx}
\usepackage{tikz}
\usepackage{MnSymbol,wasysym}
\usepackage{multicol}
\usepackage[makeroom]{cancel}
\usepackage{mathtools}
\usepackage{soul}
\usepackage{cancel}
\usepackage{url}
\usepackage{etoolbox}
\usepackage{relsize}
\usepackage{multicol}
\usepackage{caption}
\newenvironment{Figure}
  {\par\medskip\noindent\minipage{\linewidth}}
  {\endminipage\par\medskip}
\usepackage[%
    colorlinks=true,
    pdfborder={0 0 0},
    linkcolor=red
]{hyperref}
\usepackage[
    style=alphabetic,
    sorting=none 
]{biblatex}
\makeatletter
\pretocmd\start@align{
  \if@minipage\kern-\topskip\kern-\abovedisplayskip\fi
}{}{}
\makeatother

\newcommand\Ccancel[2][black]{\renewcommand\CancelColor{\color{#1}}\cancel{#2}}

\title{
 From Boxes to Polynomials: A story of generalisation}

\newtheorem*{theorem*}{Theorem}

\newcommand{\quo}[1]{`#1'}

\makeatletter
\newcommand{\proofpart}[2]{
	\par
	\addvspace{\medskipamount}
	\noindent\underline{Part #1: #2}\par\nobreak
	\addvspace{\smallskipamount}
	\@afterheading
}
\makeatother

\author{Gypsy Akhyar, Yifan Guo, and Lihexuan Yuan}
\address{School of Mathematics and Statistics, University of Melbourne}
\email{gakhyar@student.unimelb.edu.au}
\email{yifguo2@student.unimelb.edu.au}
\email{lihexuany@student.unimelb.edu.au}

\begin{document}
	\maketitle
    \begin{abstract}
        Here we will embark on a journey starting with some ostensibly inauspicious boxes. Carefully stacking them in different ways yields amazing identities. From humble beginnings at the integer version: `how many steps does it take to get from row $i$ to row $j$?' to the first upgrade: the polynomial version, before finally reaching the final upgrade: the elliptic version. Each upgrade gives a more general theorem than before. Secretly, everything is controlled by the symmetric Macdonald polynomials. Setting $q = t$ in the Macdonald polynomial takes the elliptic version of the theorem to the polynomial version. Then, letting $t$ approach $1$ reduces the polynomial version to the integer version. All the beautiful theorems and ideas come merely from stacking boxes. \\
        
    \end{abstract}
	
	\section{Introduction}
	Today, we will embark on a journey starting with some ostensibly inauspicious boxes. These boxes  will then spur a series of fascinating generalisations related to an area of great research interest --- the Macdonald polynomials.
	
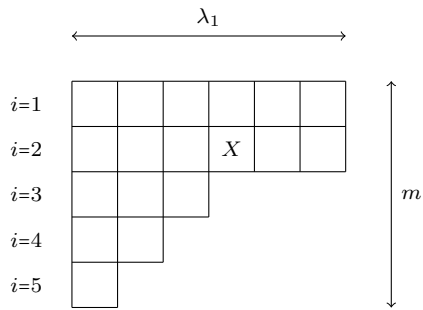
\begin{figure}[h]
         \centering
		\begin{tikzpicture}[yscale=-0.6, xscale=0.6]
			\draw (0,0) -- (6,0);
			\draw (0,1) -- (6,1);
			\draw (0,2) -- (6,2);
			\draw (0,3) -- (3,3);
			\draw (0,4) -- (2,4);
			\draw (0,5) -- (1,5);
			
			\draw (0,0) -- (0,5);
			\draw (1,0) -- (1,5);
			\draw (2,0) -- (2,4);
			\draw (3,0) -- (3,3);
			\draw (4,0) -- (4,2);
			\draw (5,0) -- (5,2);
			\draw (6,0) -- (6,2);
			\draw[<->] (7,0) -- node[anchor=west]{$\scriptstyle{m}$} (7,5) ;
			\draw[<->] (0,-1) -- node[anchor=south]{$\scriptstyle{\lambda_1}$} (6,-1) ;
			\draw (-1,0.5)  node{$\scriptstyle{i=1}$};
			\draw (-1,1.5) node{$\scriptstyle{i=2}$};
			\draw (-1,2.5) node{$\scriptstyle{i=3}$};
			\draw (-1,3.5) node{$\scriptstyle{i=4}$};
			\draw (-1,4.5) node{$\scriptstyle{i=5}$};
			
			\draw (3.5,1.5) node{$\scriptstyle{X}$};

		\end{tikzpicture}
		\caption{The box diagram for $\lambda = (6, 6, 3, 2, 1)$, with $X$ in box $(2,4)$.}
        \label{fig:Youngdiagram}
    \end{figure}

	So, what do we mean by `boxes'? The \textit{box diagram} of a sequence of non-increasing positive integers $\lambda = (\lambda_1, \lambda_2, \ldots, \lambda_m)$ is a stack of boxes with $\lambda_i$ boxes in the $i^{th}$ row \cite[Ch.\ I §1 p.\ 2]{MR1354144}. We can then index each box in the $i^{th}$ row and the $j^{th}$ column with a pair of integers $(i,j)$. Figure \ref{fig:Youngdiagram} provides such an example.

    Pick a box $b$ in $\lambda$. The \textit{content} $c(b)$ of a box is the difference between the column index $j$ and row index $i$ of the box $b$, i.e., $c(b) = j-i$ \cite[Ch.\ I §1 Ex.\ 3]{MR1354144}. An example is shown in Figure \ref{fig:contents}.
    
    \begin{figure}[h]
    \centering
			\begin{tikzpicture}[yscale=-0.5, xscale=0.5]
				\draw (0,0) -- (5,0);
				\draw (0,1) -- (5,1);
				\draw (0,2) -- (4,2);
				\draw (0,3) -- (4,3);
				\draw (0,4) -- (3,4);
				\draw (0,5) -- (2,5);
				\draw (0,0) -- (0,5);
				\draw (1,0) -- (1,5);
				\draw (2,0) -- (2,5);
				\draw (3,0) -- (3,4);
				\draw (4,0) -- (4,3);
				\draw (5,0) -- (5,1);
				
				\node at (0.5,0.5) {0};
				\node at (1.5,0.5) {1};
				\node at (2.5,0.5) {2};
				\node at (3.5,0.5) {3};
				\node at (4.5,0.5) {4};
				\node at (0.5,1.5) {-1};
				\node at (1.5,1.5) {0};
				\node at (2.5,1.5) {1};
				\node at (3.5,1.5) {2};
				\node at (0.5,2.5) {-2};
				\node at (1.5,2.5) {-1};
				\node at (2.5,2.5) {0};
				\node at (3.5,2.5) {1};
				\node at (0.5,3.5) {-3};
				\node at (1.5,3.5) {-2};
				\node at (2.5,3.5) {-1};
				\node at (0.5,4.5) {-4};
				\node at (1.5,4.5) {-3};
				
				\node at (12,2) {If $b=(5,2)$,};
				\node at (12,3) {then $c(b)=-3$.};
				
			\end{tikzpicture}
		\caption{ The contents of the boxes in $\lambda = (5, 4, 4, 3, 2)$.} \label{fig:contents} \end{figure}

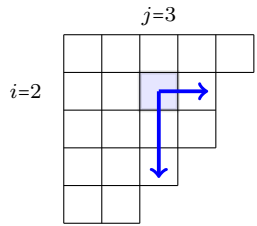
\begin{figure}[h]
    \centering
		\begin{tikzpicture}[yscale=-0.5, xscale = 0.5]
		    \filldraw [ultra thick, draw=black, fill=blue, opacity=0.1] (2,1) rectangle (3,2);
	            \draw (0,0) -- (5,0);
				\draw (0,1) -- (5,1);
				\draw (0,2) -- (4,2);
				\draw (0,3) -- (4,3);
				\draw (0,4) -- (3,4);
				\draw (0,5) -- (2,5);
				\draw (0,0) -- (0,5);
				\draw (1,0) -- (1,5);
				\draw (2,0) -- (2,5);
				\draw (3,0) -- (3,4);
				\draw (4,0) -- (4,3);
				\draw (5,0) -- (5,1);
			    \draw (-1,1.5) node{$\scriptstyle{i=2}$};
			    \draw (2.5,-0.5) node{$\scriptstyle{j=3}$};
			\draw [->, line width=0.5mm, blue] (2.5,1.5) -- (3.8,1.5);
			\draw [->, line width=0.5mm, blue]  (2.5,1.5) -- (2.5,3.8);
			
			\node at (12,1) {$b=(2,3)$};
			\node at (12,2) {\# boxes below $b = 2$};
			\node at (12,3) {\# boxes right of $b = 1$};
			\node at (12,4) {$h(b)=2+1+1=4$};
		\end{tikzpicture}
		\caption{The hook length of the box in position $(2,3)$ is 4.} \label{fig:hooklength} \end{figure}

The \textit{hook length} $h(b)$ is the number of boxes below and to the right of $b$, plus 1 for itself \cite[Ch.\ I §1 Ex.\ 1]{MR1354144}. In Figure \ref{fig:hooklength}, we can see where the term `hook length' comes from since, when drawn, the boxes counted form a hook shape.

\section{Humble beginnings: integer version.}
Let us explore quite a wonderful formula, take a peek!
\begin{theorem*}[Integer version] \emph{\cite[Ch.\ I §3 Ex.\ 4]{MR1354144}}
    \begin{equation}
		\mathlarger{\mathlarger{\prod}}_{b\in\lambda}\dfrac{n+c(b)}{h(b)}=\mathlarger{\mathlarger{\prod}}_{1\leq i<j\leq n}\dfrac{\lambda_i-\lambda_j+j-i}{j-i} \label{integer}\tag{A}
	\end{equation} 
\end{theorem*}
	
The \textbf{left-hand side} of (\ref{integer}) involves hook length and content, the notions we introduced earlier. Along with this, there is also $n$, which can be any integer greater than or equal to the length of our sequence $\lambda$.\par

We will now give an example of calculating the left-hand side for $\lambda=(5,4,4,3,2)$ with $n=5$. Since we have already shown how to calculate each piece individually, the final calculation proceeds as follows.

\begin{equation}
		\prod_{b\in\lambda}\dfrac{n+c(b)}{h(b)}=%
		\dfrac{%
			\begin{tikzpicture}[yscale=-0.45, xscale=0.45]
				\draw (0,0) -- (5,0);
				\draw (0,1) -- (5,1);
				\draw (0,2) -- (4,2);
				\draw (0,3) -- (4,3);
				\draw (0,4) -- (3,4);
				\draw (0,5) -- (2,5);
				\draw (0,0) -- (0,5);
				\draw (1,0) -- (1,5);
				\draw (2,0) -- (2,5);
				\draw (3,0) -- (3,4);
				\draw (4,0) -- (4,3);
				\draw (5,0) -- (5,1);
				\node at (0.5,0.5) {5};
				\node at (1.5,0.5) {6};
				\node at (2.5,0.5) {7};
				\node at (3.5,0.5) {8};
				\node at (4.5,0.5) {9};
				\node at (0.5,1.5) {4};
				\node at (1.5,1.5) {5};
				\node at (2.5,1.5) {6};
				\node at (3.5,1.5) {7};
				\node at (0.5,2.5) {3};
				\node at (1.5,2.5) {4};
				\node at (2.5,2.5) {5};
				\node at (3.5,2.5) {6};
				\node at (0.5,3.5) {2};
				\node at (1.5,3.5) {3};
				\node at (2.5,3.5) {4};
				\node at (0.5,4.5) {1};
				\node at (1.5,4.5) {2};
			\end{tikzpicture}\hfill
		}
		{%
			\begin{tikzpicture}[yscale=-0.45, xscale=0.45]
				\draw (0,0) -- (5,0);
				\draw (0,1) -- (5,1);
				\draw (0,2) -- (4,2);
				\draw (0,3) -- (4,3);
				\draw (0,4) -- (3,4);
				\draw (0,5) -- (2,5);
				\draw (0,0) -- (0,5);
				\draw (1,0) -- (1,5);
				\draw (2,0) -- (2,5);
				\draw (3,0) -- (3,4);
				\draw (4,0) -- (4,3);
				\draw (5,0) -- (5,1);
				\node at (0.5,0.5) {9};
				\node at (1.5,0.5) {8};
				\node at (2.5,0.5) {6};
				\node at (3.5,0.5) {4};
				\node at (4.5,0.5) {1};
				\node at (0.5,1.5) {7};
				\node at (1.5,1.5) {6};
				\node at (2.5,1.5) {4};
				\node at (3.5,1.5) {2};
				\node at (0.5,2.5) {6};
				\node at (1.5,2.5) {5};
				\node at (2.5,2.5) {3};
				\node at (3.5,2.5) {1};
				\node at (0.5,3.5) {4};
				\node at (1.5,3.5) {3};
				\node at (2.5,3.5) {1};
				\node at (0.5,4.5) {2};
				\node at (1.5,4.5) {1};
			\end{tikzpicture}\hfill
		}
		=
		\dfrac{%
			\begin{tikzpicture}[yscale=-0.45, xscale=0.45]
				\node at (0.5,0.5) {5};
				\node at (1.5,0.5) {6};
				\node at (2.5,0.5) {7};
				\node at (3.5,0.5) {8};
				\node at (4.5,0.5) {9};
				\node at (0.5,1.5) {4};
				\node at (1.5,1.5) {5};
				\node at (2.5,1.5) {6};
				\node at (3.5,1.5) {7};
				\node at (0.5,2.5) {3};
				\node at (1.5,2.5) {4};
				\node at (2.5,2.5) {5};
				\node at (3.5,2.5) {6};
				\node at (0.5,3.5) {2};
				\node at (1.5,3.5) {3};
				\node at (2.5,3.5) {4};
				\node at (0.5,4.5) {1};
				\node at (1.5,4.5) {2};
				\node at (1,0.5) {$\cdot$};
				\node at (2,0.5) {$\cdot$};
				\node at (3,0.5) {$\cdot$};
				\node at (4,0.5) {$\cdot$};
				\node at (5,0.5) {$\cdot$};
				\node at (1,1.5) {$\cdot$};
				\node at (2,1.5) {$\cdot$};
				\node at (3,1.5) {$\cdot$};
				\node at (4,1.5) {$\cdot$};
				\node at (1,2.5) {$\cdot$};
				\node at (2,2.5) {$\cdot$};
				\node at (3,2.5) {$\cdot$};
				\node at (4,2.5) {$\cdot$};
				\node at (1,3.5) {$\cdot$};
				\node at (2,3.5) {$\cdot$};
				\node at (3,3.5) {$\cdot$};
				\node at (1,4.5) {$\cdot$};
			\end{tikzpicture}\hfill
		}
		{%
			\begin{tikzpicture}[yscale=-0.45, xscale=0.45]
				\node at (0.5,0.5) {9};
				\node at (1.5,0.5) {8};
				\node at (2.5,0.5) {6};
				\node at (3.5,0.5) {4};
				\node at (4.5,0.5) {1};
				\node at (0.5,1.5) {7};
				\node at (1.5,1.5) {6};
				\node at (2.5,1.5) {4};
				\node at (3.5,1.5) {2};
				\node at (0.5,2.5) {6};
				\node at (1.5,2.5) {5};
				\node at (2.5,2.5) {3};
				\node at (3.5,2.5) {1};
				\node at (0.5,3.5) {4};
				\node at (1.5,3.5) {3};
				\node at (2.5,3.5) {1};
				\node at (0.5,4.5) {2};
				\node at (1.5,4.5) {1};
				\node at (1,0.5) {$\cdot$};
				\node at (2,0.5) {$\cdot$};
				\node at (3,0.5) {$\cdot$};
				\node at (4,0.5) {$\cdot$};
				\node at (5,0.5) {$\cdot$};
				\node at (1,1.5) {$\cdot$};
				\node at (2,1.5) {$\cdot$};
				\node at (3,1.5) {$\cdot$};
				\node at (4,1.5) {$\cdot$};
				\node at (1,2.5) {$\cdot$};
				\node at (2,2.5) {$\cdot$};
				\node at (3,2.5) {$\cdot$};
				\node at (4,2.5) {$\cdot$};
				\node at (1,3.5) {$\cdot$};
				\node at (2,3.5) {$\cdot$};
				\node at (3,3.5) {$\cdot$};
				\node at (1,4.5) {$\cdot$};
			\end{tikzpicture}\hfill
		}\label{LHS}
		\tag{B}
		\end{equation}

We now turn our attention to the \textbf{right-hand side} of (\ref{integer}). First, we want to understand how to interpret the numerator $\lambda_i - \lambda_j + j - i$.

\begin{multicols}{2}
Taking a pick of two rows, $i$ and $j$, as shown in Figure \ref{fig:rowsteps}, ask yourself `how many \textit{steps} does it take to get from row $i$ to row $j$?'.

Indeed! The answer to this question is the numerator $\lambda_i - \lambda_j + j - i$. Notice that this path requires $\lambda_i - \lambda_j$ steps to the left and $j-i$ steps down. This gives us a nice interpretation for the value $\lambda_i-\lambda_j+j-i$.

\columnbreak
	\begin{Figure} \centering
		\begin{tikzpicture}[xscale = 0.5, yscale=-0.5]
				\draw (0,0) -- (5,0);
				\draw (0,1) -- (5,1);
				\draw (0,2) -- (4,2);
				\draw (0,3) -- (4,3);
				\draw (0,4) -- (3,4);
				\draw (0,5) -- (2,5);
				\draw (0,0) -- (0,5);
				
				\draw (1,0) -- (1,5);
				\draw (2,0) -- (2,5);
				\draw (3,0) -- (3,4);
				\draw (4,0) -- (4,3);
				\draw (5,0) -- (5,1);
			\draw (-1,4.5)  node {$j=5$};
			\draw (-1,1.5) node {$i=2$};

			\draw [line width=0.5mm, blue ] (4,1) -- (4,2);
			\draw [line width=0.5mm, blue ] (2,4) -- (2,5);
			\draw (4.2,1.5)  node (a) {$\scriptstyle{0}$};
			\draw (4.2,2.5) node (b) {$\scriptstyle{1}$};
			\draw (3.5,3.2) node (c) {$\scriptstyle{2}$};
			\draw (3.2,3.6) node (d) {$\scriptstyle{3}$};
			\draw (2.5,4.2) node (e) {$\scriptstyle{4}$};
			\draw (2.2,4.6)  node (f) {$\scriptstyle{5}$};
			\draw[->] (a)  to [out=5,in=20] (b);
			\draw[->] (b)  to [out=5,in=20] (c);
			\draw[->] (c)  to [out=5,in=20] (d);
			\draw[->] (d)  to [out=5,in=20] (e);
			\draw[->] (e)  to [out=5,in=20] (f);
		\end{tikzpicture}
	
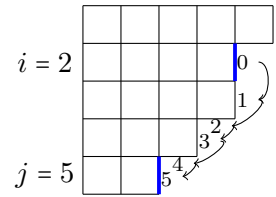
\captionof{figure}{Five steps are needed to travel from the $2^{nd}$ row to the $5^{th}$ row.} \label{fig:rowsteps} \end{Figure}
\end{multicols}

Let us now turn to the index constraints on the right-hand side of  (\ref{integer}), $1\leq i<j\leq n$. For convenience, we can use the table in Figure \ref{fig:RHSindex} to identify all eligible pairs of $(i,j)$. 
		\begin{figure}[h] \centering
			\begin{tikzpicture}[xscale = 0.5, yscale=-0.5]
				\draw[color=lightgray] (-1,-1) -- (-1,5);
				\draw[color=lightgray] (0,-1) -- (0,5);
				\draw[color=lightgray] (1,-1) -- (1,5);
				\draw[color=lightgray] (2,-1) -- (2,5);
				\draw[color=lightgray] (3,-1) -- (3,5);
				\draw[color=lightgray] (4,-1) -- (4,5);
				\draw[color=lightgray] (-2,0) -- (4,0);
				\draw[color=lightgray] (-2,1) -- (4,1);
				\draw[color=lightgray] (-2,2) -- (4,2);
				\draw[color=lightgray] (-2,3) -- (4,3);
				\draw[color=lightgray] (-2,4) -- (4,4);
				\draw[color=lightgray] (-2,5) -- (4,5);
				\draw[color=lightgray] (-2,-1) -- (-1,0);
				\draw[color=lightgray] (-1.2,-0.8)  node{$\scriptstyle{j}$};
				\draw[color=lightgray] (-1.8,-0.2)  node{$\scriptstyle{i}$};
				\draw[color=lightgray] (-0.5,-0.5)  node{$\scriptstyle{1}$};
				\draw[color=lightgray] (0.5,-0.5)  node{$\scriptstyle{2}$};
				\draw[color=lightgray] (1.5,-0.5)  node{$\scriptstyle{3}$};
				\draw[color=lightgray] (2.5,-0.5)  node{$\scriptstyle{4}$};
				\draw[color=lightgray] (3.5,-0.5)  node{$\scriptstyle{5}$};
				\draw[color=lightgray] (-1.5,0.5)  node{$\scriptstyle{1}$};
				\draw[color=lightgray] (-1.5,1.5)  node{$\scriptstyle{2}$};
				\draw[color=lightgray] (-1.5,2.5)  node{$\scriptstyle{3}$};
				\draw[color=lightgray] (-1.5,3.5)  node{$\scriptstyle{4}$};
				\draw[color=lightgray] (-1.5,4.5)  node{$\scriptstyle{5}$};
				\draw (0,0) -- (4,0);
				\draw (0,1) -- (4,1);
				\draw (1,2) -- (4,2);
				\draw (2,3) -- (4,3);
				\draw (3,4) -- (4,4);
				\draw (0,0) -- (0,1);
				\draw (1,0) -- (1,2);
				\draw (2,0) -- (2,3);
				\draw (3,0) -- (3,4);
				\draw (4,0) -- (4,4);
			\end{tikzpicture}
		\caption{A table layout which we will use for calculating the right-hand side.} \label{fig:RHSindex} \end{figure}

Returning to our example, we can finish the second half of our computation for $\lambda=(5,4,4,3,2)$.
\[
		\prod_{1\leq i<j\leq n}\dfrac{\lambda_i-\lambda_j+j-i}{j-i}=%
		\dfrac{%
			\begin{tikzpicture}[xscale = 0.5, yscale=-0.5]
				\draw (0,0) -- (4,0);
				\draw (0,1) -- (4,1);
				\draw (1,2) -- (4,2);
				\draw (2,3) -- (4,3);
				\draw (3,4) -- (4,4);
				\draw (0,0) -- (0,1);
				\draw (1,0) -- (1,2);
				\draw (2,0) -- (2,3);
				\draw (3,0) -- (3,4);
				\draw (4,0) -- (4,4);
				\draw (0.5,0.5)  node{$2$};
				\draw (1.5,0.5)  node{$3$};
				\draw (2.5,0.5)  node{$5$};
				\draw (3.5,0.5)  node{$7$};
				\draw (1.5,1.5)  node{$1$};
				\draw (2.5,1.5)  node{$3$};
				\draw (3.5,1.5)  node{$5$};
				\draw (2.5,2.5)  node{$2$};
				\draw (3.5,2.5)  node{$4$};
				\draw (3.5,3.5)  node{$2$};
			\end{tikzpicture}
		}
		{%
			\begin{tikzpicture}[xscale = 0.5, yscale=-0.5]
				\draw (0,0) -- (4,0);
				\draw (0,1) -- (4,1);
				\draw (1,2) -- (4,2);
				\draw (2,3) -- (4,3);
				\draw (3,4) -- (4,4);
				\draw (0,0) -- (0,1);
				\draw (1,0) -- (1,2);
				\draw (2,0) -- (2,3);
				\draw (3,0) -- (3,4);
				\draw (4,0) -- (4,4);
				\draw (0.5,0.5)  node{$1$};
				\draw (1.5,0.5)  node{$2$};
				\draw (2.5,0.5)  node{$3$};
				\draw (3.5,0.5)  node{$4$};
				\draw (1.5,1.5)  node{$1$};
				\draw (2.5,1.5)  node{$2$};
				\draw (3.5,1.5)  node{$3$};
				\draw (2.5,2.5)  node{$1$};
				\draw (3.5,2.5)  node{$2$};
				\draw (3.5,3.5)  node{$1$};
			\end{tikzpicture}
		}
		=%
		\dfrac{%
			\begin{tikzpicture}[xscale = 0.5, yscale=-0.5]
				\draw (0.5,0.5)  node{$2$};
				\draw (1.5,0.5)  node{$3$};
				\draw (2.5,0.5)  node{$5$};
				\draw (3.5,0.5)  node{$7$};
				\draw (1.5,1.5)  node{$1$};
				\draw (2.5,1.5)  node{$3$};
				\draw (3.5,1.5)  node{$5$};
				\draw (2.5,2.5)  node{$2$};
				\draw (3.5,2.5)  node{$4$};
				\draw (3.5,3.5)  node{$2$};
				\node at (1,0.5) {$\cdot$};
				\node at (2,0.5) {$\cdot$};
				\node at (3,0.5) {$\cdot$};
				\node at (4,0.5) {$\cdot$};
				\node at (2,1.5) {$\cdot$};
				\node at (3,1.5) {$\cdot$};
				\node at (4,1.5) {$\cdot$};
				\node at (3,2.5) {$\cdot$};
				\node at (4,2.5) {$\cdot$};
			\end{tikzpicture}\hfill
		}
		{%
			\begin{tikzpicture}[xscale = 0.5, yscale=-0.5]
				\draw (0.5,0.5)  node{$1$};
				\draw (1.5,0.5)  node{$2$};
				\draw (2.5,0.5)  node{$3$};
				\draw (3.5,0.5)  node{$4$};
				\draw (1.5,1.5)  node{$1$};
				\draw (2.5,1.5)  node{$2$};
				\draw (3.5,1.5)  node{$3$};
				\draw (2.5,2.5)  node{$1$};
				\draw (3.5,2.5)  node{$2$};
				\draw (3.5,3.5)  node{$1$};
				\node at (1,0.5) {$\cdot$};
				\node at (2,0.5) {$\cdot$};
				\node at (3,0.5) {$\cdot$};
				\node at (4,0.5) {$\cdot$};
				\node at (2,1.5) {$\cdot$};
				\node at (3,1.5) {$\cdot$};
				\node at (4,1.5) {$\cdot$};
				\node at (3,2.5) {$\cdot$};
				\node at (4,2.5) {$\cdot$};
			\end{tikzpicture}\hfill
		}
		\]

A miracle happened! This number is the same number we get for the LHS in (\ref{LHS}).
\[\prod_{\text{boxes}}\dfrac{n+c(b)}{h(b)}=\prod_{1\leq i<j\leq n}\dfrac{\lambda_i-\lambda_j+j-i}{j-i}\]
\[\dfrac{%
				\begin{tikzpicture}[blue, yscale=-0.45, xscale=0.45]
					\node at (0.5,0.5) {5};
					\node at (1.5,0.5) {6};
					\node at (2.5,0.5) {7};
					\node at (3.5,0.5) {8};
					\node at (4.5,0.5) {9};
					\node at (0.5,1.5) {4};
					\node at (1.5,1.5) {5};
					\node at (2.5,1.5) {6};
					\node at (3.5,1.5) {7};
					\node at (0.5,2.5) {3};
					\node at (1.5,2.5) {4};
					\node at (2.5,2.5) {5};
					\node at (3.5,2.5) {6};
					\node at (0.5,3.5) {2};
					\node at (1.5,3.5) {3};
					\node at (2.5,3.5) {4};
					\node at (0.5,4.5) {1};
					\node at (1.5,4.5) {2};
					\node at (1,0.5) {$\cdot$};
					\node at (2,0.5) {$\cdot$};
					\node at (3,0.5) {$\cdot$};
					\node at (4,0.5) {$\cdot$};
					\node at (5,0.5) {$\cdot$};
					\node at (1,1.5) {$\cdot$};
					\node at (2,1.5) {$\cdot$};
					\node at (3,1.5) {$\cdot$};
					\node at (4,1.5) {$\cdot$};
					\node at (1,2.5) {$\cdot$};
					\node at (2,2.5) {$\cdot$};
					\node at (3,2.5) {$\cdot$};
					\node at (4,2.5) {$\cdot$};
					\node at (1,3.5) {$\cdot$};
					\node at (2,3.5) {$\cdot$};
					\node at (3,3.5) {$\cdot$};
					\node at (1,4.5) {$\cdot$};
				\end{tikzpicture}\hfill
			}
			{%
				\begin{tikzpicture}[blue,yscale=-0.45, xscale=0.45]
					\node at (0.5,0.5) {9};
					\node at (1.5,0.5) {8};
					\node at (2.5,0.5) {6};
					\node at (3.5,0.5) {4};
					\node at (4.5,0.5) {1};
					\node at (0.5,1.5) {7};
					\node at (1.5,1.5) {6};
					\node at (2.5,1.5) {4};
					\node at (3.5,1.5) {2};
					\node at (0.5,2.5) {6};
					\node at (1.5,2.5) {5};
					\node at (2.5,2.5) {3};
					\node at (3.5,2.5) {1};
					\node at (0.5,3.5) {4};
					\node at (1.5,3.5) {3};
					\node at (2.5,3.5) {1};
					\node at (0.5,4.5) {2};
					\node at (1.5,4.5) {1};
					\node at (1,0.5) {$\cdot$};
					\node at (2,0.5) {$\cdot$};
					\node at (3,0.5) {$\cdot$};
					\node at (4,0.5) {$\cdot$};
					\node at (5,0.5) {$\cdot$};
					\node at (1,1.5) {$\cdot$};
					\node at (2,1.5) {$\cdot$};
					\node at (3,1.5) {$\cdot$};
					\node at (4,1.5) {$\cdot$};
					\node at (1,2.5) {$\cdot$};
					\node at (2,2.5) {$\cdot$};
					\node at (3,2.5) {$\cdot$};
					\node at (4,2.5) {$\cdot$};
					\node at (1,3.5) {$\cdot$};
					\node at (2,3.5) {$\cdot$};
					\node at (3,3.5) {$\cdot$};
					\node at (1,4.5) {$\cdot$};
				\end{tikzpicture}\hfill
			}
			=\dfrac{%
				\begin{tikzpicture}[red,xscale = 0.5, yscale=-0.5]
					\draw (0.5,0.5)  node{$2$};
					\draw (1.5,0.5)  node{$3$};
					\draw (2.5,0.5)  node{$5$};
					\draw (3.5,0.5)  node{$7$};
					\draw (1.5,1.5)  node{$1$};
					\draw (2.5,1.5)  node{$3$};
					\draw (3.5,1.5)  node{$5$};
					\draw (2.5,2.5)  node{$2$};
					\draw (3.5,2.5)  node{$4$};
					\draw (3.5,3.5)  node{$2$};
					\node at (1,0.5) {$\cdot$};
					\node at (2,0.5) {$\cdot$};
					\node at (3,0.5) {$\cdot$};
					\node at (4,0.5) {$\cdot$};
					\node at (2,1.5) {$\cdot$};
					\node at (3,1.5) {$\cdot$};
					\node at (4,1.5) {$\cdot$};
					\node at (3,2.5) {$\cdot$};
					\node at (4,2.5) {$\cdot$};
				\end{tikzpicture}\hfill
			}
			{%
				\begin{tikzpicture}[red, xscale = 0.5, yscale=-0.5]
					\draw (0.5,0.5)  node{$1$};
					\draw (1.5,0.5)  node{$2$};
					\draw (2.5,0.5)  node{$3$};
					\draw (3.5,0.5)  node{$4$};
					\draw (1.5,1.5)  node{$1$};
					\draw (2.5,1.5)  node{$2$};
					\draw (3.5,1.5)  node{$3$};
					\draw (2.5,2.5)  node{$1$};
					\draw (3.5,2.5)  node{$2$};
					\draw (3.5,3.5)  node{$1$};
					\node at (1,0.5) {$\cdot$};
					\node at (2,0.5) {$\cdot$};
					\node at (3,0.5) {$\cdot$};
					\node at (4,0.5) {$\cdot$};
					\node at (2,1.5) {$\cdot$};
					\node at (3,1.5) {$\cdot$};
					\node at (4,1.5) {$\cdot$};
					\node at (3,2.5) {$\cdot$};
					\node at (4,2.5) {$\cdot$};
				\end{tikzpicture}\hfill
			}\]
			We can see that the two sides are equal by comparing the numbers in the cross-multiplication row by row.
			\begin{equation}
			\begin{tikzpicture}[blue, yscale=-0.45, xscale=0.45]
				\node at (0.5,0.5) {5};
				\node at (1.5,0.5) {6};
				\node at (2.5,0.5) {7};
				\node at (3.5,0.5) {8};
				\node at (4.5,0.5) {9};
				\node at (0.5,1.5) {4};
				\node at (1.5,1.5) {5};
				\node at (2.5,1.5) {6};
				\node at (3.5,1.5) {7};
				\node at (0.5,2.5) {3};
				\node at (1.5,2.5) {4};
				\node at (2.5,2.5) {5};
				\node at (3.5,2.5) {6};
				\node at (0.5,3.5) {2};
				\node at (1.5,3.5) {3};
				\node at (2.5,3.5) {4};
				\node at (0.5,4.5) {1};
				\node at (1.5,4.5) {2};
				\node at (1,0.5) {$\cdot$};
				\node at (2,0.5) {$\cdot$};
				\node at (3,0.5) {$\cdot$};
				\node at (4,0.5) {$\cdot$};
				\node at (5,0.5) {$\cdot$};
				\node at (1,1.5) {$\cdot$};
				\node at (2,1.5) {$\cdot$};
				\node at (3,1.5) {$\cdot$};
				\node at (4,1.5) {$\cdot$};
				\node at (1,2.5) {$\cdot$};
				\node at (2,2.5) {$\cdot$};
				\node at (3,2.5) {$\cdot$};
				\node at (4,2.5) {$\cdot$};
				\node at (1,3.5) {$\cdot$};
				\node at (2,3.5) {$\cdot$};
				\node at (3,3.5) {$\cdot$};
				\node at (1,4.5) {$\cdot$};
			\end{tikzpicture}\times
			\begin{tikzpicture}[red, xscale = 0.5, yscale=-0.5]
				\draw (0.5,0.5)  node{$1$};
				\draw (1.5,0.5)  node{$2$};
				\draw (2.5,0.5)  node{$3$};
				\draw (3.5,0.5)  node{$4$};
				\draw (1.5,1.5)  node{$1$};
				\draw (2.5,1.5)  node{$2$};
				\draw (3.5,1.5)  node{$3$};
				\draw (2.5,2.5)  node{$1$};
				\draw (3.5,2.5)  node{$2$};
				\draw (3.5,3.5)  node{$1$};
				\filldraw [opacity=0] (3,4) rectangle (3.3,4.6);
				\node at (1,0.5) {$\cdot$};
				\node at (2,0.5) {$\cdot$};
				\node at (3,0.5) {$\cdot$};
				\node at (4,0.5) {$\cdot$};
				\node at (2,1.5) {$\cdot$};
				\node at (3,1.5) {$\cdot$};
				\node at (4,1.5) {$\cdot$};
				\node at (3,2.5) {$\cdot$};
				\node at (4,2.5) {$\cdot$};
			\end{tikzpicture}
			=
			\begin{tikzpicture}[blue,yscale=-0.45, xscale=0.45]
				\node at (0.5,0.5) {9};
				\node at (1.5,0.5) {8};
				\node at (2.5,0.5) {6};
				\node at (3.5,0.5) {4};
				\node at (4.5,0.5) {1};
				\node at (0.5,1.5) {7};
				\node at (1.5,1.5) {6};
				\node at (2.5,1.5) {4};
				\node at (3.5,1.5) {2};
				\node at (0.5,2.5) {6};
				\node at (1.5,2.5) {5};
				\node at (2.5,2.5) {3};
				\node at (3.5,2.5) {1};
				\node at (0.5,3.5) {4};
				\node at (1.5,3.5) {3};
				\node at (2.5,3.5) {1};
				\node at (0.5,4.5) {2};
				\node at (1.5,4.5) {1};
				\node at (1,0.5) {$\cdot$};
				\node at (2,0.5) {$\cdot$};
				\node at (3,0.5) {$\cdot$};
				\node at (4,0.5) {$\cdot$};
				\node at (5,0.5) {$\cdot$};
				\node at (1,1.5) {$\cdot$};
				\node at (2,1.5) {$\cdot$};
				\node at (3,1.5) {$\cdot$};
				\node at (4,1.5) {$\cdot$};
				\node at (1,2.5) {$\cdot$};
				\node at (2,2.5) {$\cdot$};
				\node at (3,2.5) {$\cdot$};
				\node at (4,2.5) {$\cdot$};
				\node at (1,3.5) {$\cdot$};
				\node at (2,3.5) {$\cdot$};
				\node at (3,3.5) {$\cdot$};
				\node at (1,4.5) {$\cdot$};
			\end{tikzpicture}\times
			\begin{tikzpicture}[red, xscale = 0.5, yscale=-0.5]
				\draw (0.5,0.5)  node{$2$};
				\draw (1.5,0.5)  node{$3$};
				\draw (2.5,0.5)  node{$5$};
				\draw (3.5,0.5)  node{$7$};
				\draw (1.5,1.5)  node{$1$};
				\draw (2.5,1.5)  node{$3$};
				\draw (3.5,1.5)  node{$5$};
				\draw (2.5,2.5)  node{$2$};
				\draw (3.5,2.5)  node{$4$};
				\draw (3.5,3.5)  node{$2$};
				\filldraw [opacity=0] (3,4) rectangle (3.3,4.6);
				\node at (1,0.5) {$\cdot$};
				\node at (2,0.5) {$\cdot$};
				\node at (3,0.5) {$\cdot$};
				\node at (4,0.5) {$\cdot$};
				\node at (2,1.5) {$\cdot$};
				\node at (3,1.5) {$\cdot$};
				\node at (4,1.5) {$\cdot$};
				\node at (3,2.5) {$\cdot$};
				\node at (4,2.5) {$\cdot$};
			\end{tikzpicture}\label{cross-multiply}
			\tag{C}
			\end{equation} 

Each row forms a factorial! Why is this happening? Does it always work?

\begin{figure}[h] \centering
		\begin{tikzpicture}[yscale=-1]
		    \filldraw [fill=blue, opacity=0.1] (1,0) rectangle (2,1);
		    \draw[ultra thick] (1.5,0.5) -- (4.8,0.5);
		    \draw[ultra thick] (1.5,0.5) -- (1.5,4.8);
		    \draw[ultra thick] (5,1) -- (4,1) -- (4,2) -- (4,3) -- (3,3)-- (3,4) -- (2,4) --(2,5) -- (1,5);
			\draw (0,0) -- (5,0);
			\draw (0,1) -- (5,1);
			\draw (0,2) -- (4,2);
			\draw (0,3) -- (4,3);
			\draw (0,4) -- (3,4);
			\draw (0,5) -- (2,5);
			\draw (0,0) -- (0,5);
			\draw (1,0) -- (1,5);
			\draw (2,0) -- (2,5);
			\draw (3,0) -- (3,4);
			\draw (4,0) -- (4,3);
			\draw (5,0) -- (5,1);
			\draw (5.2,0.5)  node (a) {$\scriptstyle{0}$};
			\draw (4.5,1.2) node (b) {$\scriptstyle{1}$};
			\draw (4.2,1.5) node (c) {$\scriptstyle{2}$};
			\draw (4.2,2.5) node (d) {$\scriptstyle{3}$};
			\draw (3.6,3.2)  node (e) {$\scriptstyle{4}$};
			\draw (3.2,3.5)  node (f) {$\scriptstyle{5}$};
			\draw (2.5,4.2)  node (g) {$\scriptstyle{6}$};
			\draw (2.2,4.5)  node (h) {$\scriptstyle{7}$};
			\draw (1.5,5.2)  node (i) {$\scriptstyle{8}$};
			\draw[->] (a)  to [out=5,in=20] (b);
			\draw[->] (b)  to [out=5,in=20] (c);
			\draw[->] (c)  to [out=5,in=20] (d);
			\draw[->] (d)  to [out=5,in=20] (e);
			\draw[->] (e)  to [out=5,in=20] (f);
			\draw[->] (f)  to [out=5,in=20] (g);
			\draw[->] (g)  to [out=5,in=20] (h);
			\draw[->] (h)  to [out=5,in=20] (i);
		\end{tikzpicture}
		\caption{The hook length of the box at $(1,2)$ equals the corresponding path length along the boundary.} \label{fig:boundary} \end{figure}
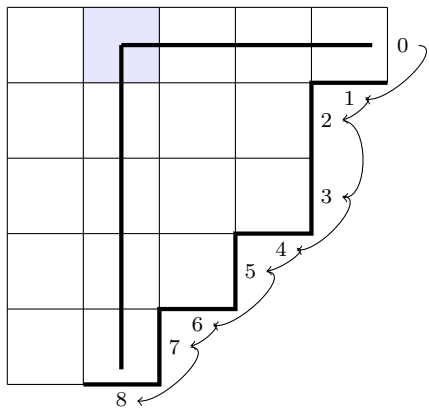

Let us now explain the key intuition behind the proof of this formula. To do this, we first interpret the hook length in a slightly different way. Instead of counting the number of boxes covered by the inverted `L' shape, we travel along the boundary, as shown in Figure \ref{fig:boundary}.

\pagebreak
For any box $(i,j)$, its hook length answers the question:
\begin{center}
    \quo{\textit{How many steps does it take to travel from the right end of the $i^{\text{th}}$ row to the \textbf{bottom} of the $j^{\text{th}}$ column?}}
\end{center}

This feels suspiciously similar to our interpretation of the numerator, $\lambda_i-\lambda_j+j-i$, which recall was the answer to the question:
\begin{center}
    \quo{\textit{How many steps does it take to travel from the $i^{\text{th}}$ row to the (\textbf{right end of the}) $j^{\text{th}}$ row?}}
\end{center}

So, perhaps $h(b)$ and $\lambda_i-\lambda_j+j-i$ are based on a more general notion? In hopes of finding this, we take a deeper look at what they describe. Now, fix some row $i$.

\begin{multicols}{2}
	Iterating the hook length over the boxes
	\[\dprod_{\text{boxes in row $i$}} h(b)\] describes the number of steps to travel from the $i^{th}$ row to the bottom of each column, as shown in Figure \ref{fig:bottomcol}.
	\begin{Figure} \centering
		\begin{tikzpicture}[xscale = 0.7, yscale=-0.7]
			\draw (0,0) -- (5,0);
			\draw (0,1) -- (5,1);
			\draw (0,2) -- (4,2);
			\draw (0,3) -- (4,3);
			\draw (0,4) -- (3,4);
			\draw (0,5) -- (2,5);
			\draw (0,0) -- (0,5);
			\draw (1,0) -- (1,5);
			\draw (2,0) -- (2,5);
			\draw (3,0) -- (3,4);
			\draw (4,0) -- (4,3);
			\draw (5,0) -- (5,1);
			\draw (3.6,3.2)  node (e) {$\scriptstyle{2}$};
			\draw (2.5,4.2)  node (g) {$\scriptstyle{4}$};
			\draw (1.5,5.2)  node (i) {$\scriptstyle{6}$};
			\draw (0.5,5.2)  node (i) {$\scriptstyle{7}$};
			\draw (4.6,1.5) node {start};
			\draw [line width=0.5mm, blue ] (3,3) -- (4,3);
			\draw [line width=0.5mm, blue ] (2,4) -- (3,4);
			\draw [line width=0.5mm, blue ] (0,5) -- (2,5);
		\end{tikzpicture}
	\captionof{figure}{Starting at $i=2$ and travelling to each column.} \label{fig:bottomcol} \end{Figure}
\columnbreak
	Iterating $\lambda_i-\lambda_j+j-i$ such that $1\leq i<j \leq n$
\[\dprod_{1\leq i<j \leq n} \lambda_i-\lambda_j+j-i\] describes the number of steps to travel from the $i^{th}$ row to the right end of each subsequent row, as shown in Figure \ref{fig:rightrow}.
\begin{Figure} \centering
	\begin{tikzpicture}[xscale = 0.7, yscale=-0.7]
		\draw (0,0) -- (5,0);
		\draw (0,1) -- (5,1);
		\draw (0,2) -- (4,2);
		\draw (0,3) -- (4,3);
		\draw (0,4) -- (3,4);
		\draw (0,5) -- (2,5);
		\draw (0,0) -- (0,5);
		\draw (1,0) -- (1,5);
		\draw (2,0) -- (2,5);
		\draw (3,0) -- (3,4);
		\draw (4,0) -- (4,3);
		\draw (5,0) -- (5,1);
		\draw (4.2,2.5) node (d) {$\scriptstyle{1}$};
		\draw (3.2,3.5)  node (f) {$\scriptstyle{3}$};
		\draw (2.2,4.5)  node (h) {$\scriptstyle{5}$};
		\draw (4.6,1.5) node {start};
		\draw [line width=0.5mm, blue ] (4,2) -- (4,3);
		\draw [line width=0.5mm, blue ] (3,3) -- (3,4);
		\draw [line width=0.5mm, blue ] (2,4) -- (2,5);
	\end{tikzpicture}
\captionof{figure}{Starting at $i=2$ and travelling to each row.} \label{fig:rightrow} \end{Figure}
\end{multicols}

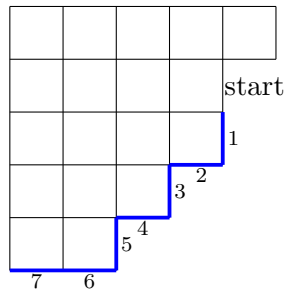
\begin{figure}[h]
    \centering
    \begin{center}
		\begin{tikzpicture}[xscale = 0.7, yscale=-0.7]
    		\draw (0,0) -- (5,0);
    		\draw (0,1) -- (5,1);
    		\draw (0,2) -- (4,2);
    		\draw (0,3) -- (4,3);
    		\draw (0,4) -- (3,4);
    		\draw (0,5) -- (2,5);
    		\draw (0,0) -- (0,5);
    		\draw (1,0) -- (1,5);
    		\draw (2,0) -- (2,5);
    		\draw (3,0) -- (3,4);
    		\draw (4,0) -- (4,3);
    		\draw (5,0) -- (5,1);
    		\draw (4.2,2.5) node (d) {$\scriptstyle{1}$};
    		\draw (3.2,3.5)  node (f) {$\scriptstyle{3}$};
    		\draw (2.2,4.5)  node (h) {$\scriptstyle{5}$};
    		\draw (3.6,3.2)  node (e) {$\scriptstyle{2}$};
			\draw (2.5,4.2)  node (g) {$\scriptstyle{4}$};
			\draw (1.5,5.2)  node (i) {$\scriptstyle{6}$};
			\draw (0.5,5.2)  node (i) {$\scriptstyle{7}$};
			\draw (4.6,1.5) node {start};
    		\draw [line width=0.5mm, blue ] (4,2) -- (4,3);
    		\draw [line width=0.5mm, blue ] (3,3) -- (3,4);
    		\draw [line width=0.5mm, blue ] (2,4) -- (2,5);
			\draw [line width=0.5mm, blue ] (3,3) -- (4,3);
			\draw [line width=0.5mm, blue ] (2,4) -- (3,4);
			\draw [line width=0.5mm, blue ] (0,5) -- (2,5);
		\end{tikzpicture}
	\caption{Combining both sides gives a complete path along the boundary.}
	\label{combinepath}
	\end{center}
\end{figure}

We can now see that in Figure \ref{combinepath}, these will always give us the consecutive integers $1, 2, 3,\dots$ as we travel along the complete boundary. That is where these factorials come from!

	Thus we obtain the theorem inappro (\ref{integer}). Things are going to get a whole lot more interesting!
	
\section{The first upgrade: polynomial version.}
The key step in our proof for the integer version is that, after cross-multiplication, both sides consisted of the same factors in equation (\ref{cross-multiply}). If we were to replace each number $k$ with a polynomial $1-t^k$, both sides should still contain the same polynomial factors. For example, when $2\cdot5\cdot6 = 6\cdot2\cdot5$ gets replaced by $ (1-t^2)(1-t^5)(1-t^6) = (1-t^6)(1-t^2)(1-t^5)$, it is still true. This is the idea for our first generalisation.

So, the polynomial analogue of equation (\ref{cross-multiply}) gives us:

\begin{theorem*}[Polynomial version I] \emph{\cite[Ch.\ I §3 Ex.\ 1]{MR1354144}} 
			\begin{equation}
				\prod_{b\in \lambda} \dfrac{1-t^{n+c(b)}}{1-t^{h(b)}} = \prod_{1\leq i <j \leq n} \dfrac{1-t^{\lambda_i-\lambda_j+j-i}}{1-t^{j-i}}.
				\label{polynomial}\tag{D}
			\end{equation}
			
\end{theorem*} 

Observe that $\frac{1-t^k}{1-t} = 1+t+t^2+\ldots+t^{k-1}$ approaches $k$ as $t \rightarrow 1$, and hence the polynomial version will collapse to the integer version of the formula!
		
To prepare for our next upgrade, we first take a closer look at the content $c(b)$ and hook length $h(b)$ of a box in equation (\ref{polynomial}). Let us introduce some new terminology.

\begin{multicols}{2}
    \begin{Figure}
        \centering
        \begin{tikzpicture}[yscale=-1]
				\draw (0,0) -- (5,0) -- (5,1) -- (4,1) -- (4,3)  -- (2,3) -- (2,4) -- (1,4) -- (1,5) -- (0,5) -- (0,0);
				\draw[<->] (2.5,0) -- node[anchor=west]{$\scriptstyle{coleg_\lambda(b)}$} (2.5,1.3) ;
				\draw[<->] (2.5,1.7) -- node[anchor=west]{$\scriptstyle{leg_\lambda(b)}$} (2.5,3) ;
				\draw[<->] (0,1.5) -- node[anchor=south]{$\scriptstyle{coarm_\lambda(b)}$} (2.3,1.5) ;
				\draw[<->] (2.7,1.5) -- node[anchor=south]{$\scriptstyle{arm_\lambda(b)}$} (4,1.5) ;
				\draw (2.5,1.5)  node[shape=rectangle,draw]{$b$};
		\end{tikzpicture}
        
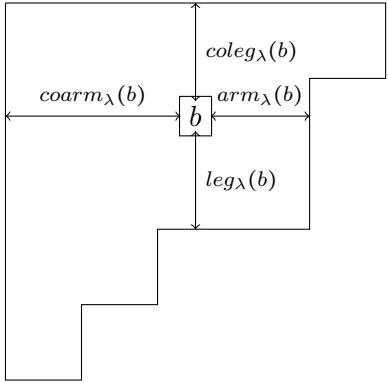
\captionof{figure}{The coarm, coleg, arm, and leg of a box.}
        \label{fig:armleg}
    \end{Figure}

		\columnbreak
			
			As indicated in Figure \ref{fig:armleg}, for a box $b$, we define the \textit{arm} and \textit{coarm} of a box to be the number of boxes to the right and to the left of $b$, respectively. 
			
			Similarly, we call the number of boxes below and above $b$ the \textit{leg} and \textit{coleg}, respectively.
			
			So, we can express the content and the hook length as: 
			\begin{align*}
				c(b) &= coarm_\lambda(b) - coleg_\lambda(b)\\
				h(b) &= arm_\lambda(b)+leg_\lambda(b)+1
			\end{align*}
\end{multicols}

		We can now rewrite the left-hand side of our theorem:
		\begin{theorem*}[Polynomial version II]
			\begin{equation*}
				\prod_{b\in \lambda} \dfrac{1-t^{n+coarm_\lambda(b)- coleg_\lambda(b)}}{1-t^{arm_\lambda(b)+leg_\lambda(b)+1}} = \prod_{1\leq i <j \leq n} \dfrac{1-t^{\lambda_i-\lambda_j+j-i}}{1-t^{j-i}} 
			\end{equation*}
		\end{theorem*}
		
		\pagebreak
				\section{The second upgrade: elliptic version.}
		Now that we have introduced our new terminology, we generalise the polynomial version with the introduction of a new variable $q$. 
		
			\begin{theorem*}[Elliptic version] \emph{\cite[Ch.\ VI (6.11$'$)]{MR1354144}}
			\begin{equation}
				\prod_{b\in \lambda} \dfrac{1-q^{coarm_\lambda(b)}t^{n-coleg_\lambda(b)}}{1-q^{arm_\lambda(b)}t^{leg_\lambda(b)+1}} 
				=
				\prod_{1\leq i < j \leq n} \prod_{r=0}^{\lambda_i-\lambda_j-1} \dfrac{1-q^rt^{j-i+1}}{1-q^rt^{j-i}}. \label{elliptic}\tag{E}
			\end{equation}
		\end{theorem*}
		
		In the special case where $q=t$, the elliptic version of the theorem (\ref{elliptic}) reduces to the polynomial version of the theorem (\ref{polynomial}). This is easy to see for the left-hand side, but the right-hand side requires some more cancellations to return to the polynomial version.
		
		\begin{equation*}
    		\begin{split}
    			\text{RHS of (\ref{elliptic}) } \rvert_{q=t} = & \left. \prod_{1\leq i < j \leq n} \prod_{r=0}^{\lambda_i-\lambda_j-1} \dfrac{1-q^rt^{j-i+1}}{1-q^rt^{j-i}} \, \right|_{q=t}\\
    			= & \prod_{1\leq i < j \leq n} \dfrac{\Ccancel[gray]{1-t^0t^{j-i+1}}}{1-t^0t^{j-i}} \dfrac{\Ccancel[gray]{1-t^1t^{j-i+1}}}{\Ccancel[gray]{1-t^1t^{j-i}}} \dots \dfrac{1-t^{\lambda_i-\lambda_j-1}t^{j-i+1}}{\Ccancel[gray]{1-t^{\lambda_i-\lambda_j-1}t^{j-i}}}\\
    			= & \prod_{1\leq i <j \leq n} \dfrac{1-t^{\lambda_i-\lambda_j+j-i}}{1-t^{j-i}}.
    		\end{split}
		\end{equation*}
		
	    Returning to our more generalised case where $q$ and $t$ are not necessarily equal, let us also prove it with a visual and combinatorial approach.

		First, observe the right-hand side of (\ref{elliptic}). Notice that the powers of $t$ in the numerator and the denominator differ by 1: \(\dfrac{1-q^rt^{j-i+1}}{1-q^rt^{j-i}}\). We will try to utilise this observation in hopes of finding a nice cancellation. Thus, our first step here will be to change the order of the indices $r$, $i$, and $j$, so that we can first take the product over $j$, instead of $r$.
		
		For convenience, let us temporarily denote this fraction as \[rK+\alpha_{ij} := \dfrac{1-q^rt^{j-i+1}}{1-q^rt^{j-i}}.\]
		Here, $K$ and $\alpha$ are not variables but are only symbols for this new notation. We can now write the right hand side of (\ref{elliptic}) as:
		\begin{equation*}
		    \prod_{1\leq i < j \leq n} \prod_{r=0}^{\lambda_i-\lambda_j-1} rK+\alpha_{ij}.
		\end{equation*}
		
		To better understand what is happening in this expression, let us return to our example $\lambda = (5, 4, 4, 3, 2)$ with a choice of $n=5$.
		
		When $i=1$ and $j=5$, we have \[0 \leq r \leq \lambda_1-\lambda_5-1 = 5-2-1 = 2.\]
		This gives us three factors, $2K + \alpha_{15}$, $K+\alpha_{15}$ and $\alpha_{15}$. We will write these three factors in the first line of the following box diagram, as shown in bold in Figure \ref{fig:RHS}. 
		
	\begin{figure}
	    \centering
	    	
		\begin{center}
			\begin{tikzpicture}[yscale=-1.2, xscale=1.2]
				\draw (0,0) -- (5,0);
				\draw (0,1) -- (5,1);
				\draw (0,2) -- (4,2);
				\draw (0,3) -- (4,3);
				\draw (0,4) -- (3,4);
				\draw (0,5) -- (2,5);
				\draw (0,0) -- (0,5);
				\draw (1,0) -- (1,5);
				\draw (2,0) -- (2,5);
				\draw (3,0) -- (3,4);
				\draw (4,0) -- (4,3);
				\draw (5,0) -- (5,1);
				\draw (-1,0.5)  node{$\scriptstyle{i=1}$};
				\draw (-1,1.5) node{$\scriptstyle{i=2}$};
				\draw (-1,2.5) node{$\scriptstyle{i=3}$};
				\draw (-1,3.5) node{$\scriptstyle{i=4}$};
				\draw (-1,4.5) node{$\scriptstyle{i=5}$};
				\draw (2.5,0.2)  node{$\boldsymbol{\scriptstyle{2K + \alpha_{15}}}$};
				\draw (3.5,0.2)  node{$\scriptstyle{\boldsymbol{K+\alpha_{15}}}$};
				\draw (4.5,0.2)  node{$\scriptstyle{\boldsymbol{\alpha_{15}}}$};
				\draw (3.5,0.4)  node{$\scriptstyle{K+\alpha_{14}}$};
				\draw (4.5,0.4)  node{$\scriptstyle{\alpha_{14}}$};
				\draw (4.5,0.6)  node{$\scriptstyle{\alpha_{13}}$};
				\draw (4.5,0.8)  node{$\scriptstyle{\alpha_{12}}$};
				\draw (2.5,1.2)  node{$\scriptstyle{K + \alpha_{25}}$};
				\draw (3.5,1.2)  node{$\scriptstyle{\alpha_{25}}$};
				\draw (3.5,1.4)  node{$\scriptstyle{\alpha_{24}}$};
				\draw (2.5,2.2)  node{$\scriptstyle{K + \alpha_{35}}$};
				\draw (3.5,2.2)  node{$\scriptstyle{\alpha_{35}}$};
				\draw (3.5,2.4)  node{$\scriptstyle{\alpha_{34}}$};
				\draw (2.5,3.2)  node{$\scriptstyle{\alpha_{45}}$};
			\end{tikzpicture}
		\end{center}
		\caption{Factors on the right hand side of the equation for $\lambda = (5, 4, 4, 3, 2)$.}
	    \label{fig:RHS}
	\end{figure}
			The diagram can help us better visualise the change of product order that we had mentioned earlier. Repeating for all allowed values of $i$, $j$ and $r$, we obtain Figure \ref{fig:RHS}.
		
		Notice how factors in each cell now have the same values of $i$ and $r$; this enables us to simplify the factors. For example, if we focus on the rightmost cell in the first row we see that we have a fair amount of cancellations:
		\begin{align*}				    \vcenter{\hbox{\begin{tikzpicture}[yscale=-1.2, xscale=1.2]
							\draw (4.5,0.2)  node{$\scriptstyle{\alpha_{15}}$};
							\draw (4.5,0.4)  node{$\scriptstyle{\alpha_{14}}$};
							\draw (4.5,0.6)  node{$\scriptstyle{\alpha_{13}}$};
							\draw (4.5,0.8)  node{$\scriptstyle{\alpha_{12}}$};
							\draw[draw=black] (4,0) rectangle ++(1,1);
				\end{tikzpicture}}}
				&=\frac{1-q^0t^5}{\Ccancel[gray]{1-q^0t^4}}
				\frac{\Ccancel[gray]{1-q^0t^4}}{\Ccancel[gray]{1-q^0t^4}}
				\frac{\Ccancel[gray]{1-q^0t^3}}{\Ccancel[gray]{1-q^0t^2}}
				\frac{\Ccancel[gray]{1-q^0t^2}}{1-q^0t^1} =
				\frac{1-q^0t^5}{1-q^0t^1}.
			\end{align*}
		 We will repeat the same cancellation process for each cell until we have the following new table of factors in Figure \ref{fig:RHScancel}. 
\begin{figure}[h]
    \centering
    \begin{center}
			\begin{tikzpicture}[yscale=-1.2, xscale=1.2]
				\draw (0,0) -- (5,0);
				\draw (0,1) -- (5,1);
				\draw (0,2) -- (4,2);
				\draw (0,3) -- (4,3);
				\draw (0,4) -- (3,4);
				\draw (0,5) -- (2,5);
				\draw (0,0) -- (0,5);
				\draw (1,0) -- (1,5);
				\draw (2,0) -- (2,5);
				\draw (3,0) -- (3,4);
				\draw (4,0) -- (4,3);
				\draw (5,0) -- (5,1);
				
				\draw (2.5,0.5)  node{$\scriptstyle{\frac{1-q^2t^5}{1-q^2t^4}}$};
				\draw (3.5,0.5)  node{$\scriptstyle{\frac{1-q^1t^5}{1-q^1t^3}}$};
				\draw (4.5,0.5)  node{$\scriptstyle{\frac{1-q^0t^5}{1-q^0t^1}}$};
				\draw (2.5,1.5)  node{$\scriptstyle{\frac{1-q^1t^4}{1-q^1t^3}}$};
				\draw (3.5,1.5)  node{$\scriptstyle{\frac{1-q^0t^4}{1-q^0t^2}}$};
				\draw (2.5,2.5)  node{$\scriptstyle{\frac{1-q^1t^3}{1-q^1t^2}}$};
				\draw (3.5,2.5)  node{$\scriptstyle{\frac{1-q^0t^3}{1-q^0t^1}}$};
				\draw (2.5,3.5)  node{$\scriptstyle{\frac{1-q^0t^2}{1-q^0t^1}}$};
			\end{tikzpicture}
		\end{center}
    \caption{The right-hand side of the equation after cancellations.}
    \label{fig:RHScancel}
\end{figure}
    \pagebreak
		Next, let us simply `play around' with the factors a little bit. 
		
		Reverse the order of all the numerators in each row as shown in Figure \ref{fig:reverseorder}.
\begin{figure}[h]
    \centering
		\begin{center}
				\begin{tikzpicture}[yscale=-1.2, xscale=1.2]
					\draw (0,0) -- (5,0);
					\draw (0,1) -- (5,1);
					\draw (0,2) -- (4,2);
					\draw (0,3) -- (4,3);
					\draw (0,4) -- (3,4);
					\draw (0,5) -- (2,5);
					\draw (0,0) -- (0,5);
					\draw (1,0) -- (1,5);
					\draw (2,0) -- (2,5);
					\draw (3,0) -- (3,4);
					\draw (4,0) -- (4,3);
					\draw (5,0) -- (5,1);

					\draw (0.5,0.5)  node{$\scriptstyle{\frac{1-q^0t^5}{1}}$};
					\draw (1.5,0.5)  node{$\scriptstyle{\frac{1-q^1t^5}{1}}$};
					\draw (2.5,0.5)  node{$\scriptstyle{\frac{1-q^2t^5}{1-q^2t^4}}$};
					\draw (3.5,0.5)  node{$\scriptstyle{\frac{1}{1-q^1t^3}}$};
					\draw (4.5,0.5)  node{$\scriptstyle{\frac{1}{1-q^0t^1}}$};
					\draw (0.5,1.5)  node{$\scriptstyle{\frac{1-q^0t^4}{1}}$};
					\draw (1.5,1.5)  node{$\scriptstyle{\frac{1-q^1t^4}{1}}$};
					\draw (2.5,1.5)  node{$\scriptstyle{\frac{1}{1-q^1t^3}}$};
					\draw (3.5,1.5)  node{$\scriptstyle{\frac{1}{1-q^0t^2}}$};
					\draw (0.5,2.5)  node{$\scriptstyle{\frac{1-q^0t^3}{1}}$};
					\draw (1.5,2.5)  node{$\scriptstyle{\frac{1-q^1t^3}{1}}$};
					\draw (2.5,2.5)  node{$\scriptstyle{\frac{1}{1-q^1t^2}}$};
					\draw (3.5,2.5)  node{$\scriptstyle{\frac{1}{1-q^0t^1}}$};
					\draw (0.5,3.5)  node{$\scriptstyle{\frac{1-q^0t^2}{1}}$};
					\draw (2.5,3.5)  node{$\scriptstyle{\frac{1}{1-q^0t^1}}$};
					
				\end{tikzpicture}
			\end{center}
    \caption{Reversing the order of numerators.}
    \label{fig:reverseorder}
\end{figure}
		Do you notice any connection between this resultant table and the left-hand side of the equation (\ref{elliptic})?
		\[\text{LHS of (\ref{elliptic})} = \prod_{b\in \lambda} \dfrac{1-q^{coarm_\lambda(b)}t^{n-coleg_\lambda(b)}}{1-q^{arm_\lambda(b)}t^{leg_\lambda(b)+1}}. \]
		
		Like the left-hand side, numerators in this table seem to also take the form $1-q^{coarm_\lambda(b)}t^{n-coleg_\lambda(b)}$. The same can also be said about the denominators which take the form $1-q^{arm_\lambda(b)}t^{leg_\lambda(b)+1}$.
		
		If only we had all the missing numerators and denominators filled in, then we would have exactly the same factors as the left-hand side and could show that the two sides of (\ref{elliptic}) are equal. We will now add the missing terms as indicated in Figure \ref{fig:addmissingfactors}.
\begin{figure}[h]
    \centering
		\begin{center}
			\begin{tikzpicture}[yscale=-1.2, xscale=1.2]
				\draw (0,0) -- (5,0);
				\draw (0,1) -- (5,1);
				\draw (0,2) -- (4,2);
				\draw (0,3) -- (4,3);
				\draw (0,4) -- (3,4);
				\draw (0,5) -- (2,5);
				\draw (0,0) -- (0,5);
				\draw (1,0) -- (1,5);
				\draw (2,0) -- (2,5);
				\draw (3,0) -- (3,4);
				\draw (4,0) -- (4,3);
				\draw (5,0) -- (5,1);

				\draw (0.5,0.5)  node{$\scriptstyle{\frac{1-q^0t^5}{\textcolor{red}{1-q^4t^5}}}$};
				\draw (1.5,0.5)  node{$\scriptstyle{\frac{1-q^1t^5}{\textcolor{red}{1-q^3t^5}}}$};
				\draw (2.5,0.5)  node{$\scriptstyle{\frac{1-q^2t^5}{1-q^2t^4}}$};
				\draw (3.5,0.5)  node{$\scriptstyle{\frac{\textcolor{red}{1-q^3t^5}}{1-q^1t^3}}$};
				\draw (4.5,0.5)  node{$\scriptstyle{\frac{\textcolor{red}{1-q^4t^5}}{1-q^0t^1}}$};
				\draw (0.5,1.5)  node{$\scriptstyle{\frac{1-q^0t^4}{\textcolor{red}{1-q^3t^4}}}$};
				\draw (1.5,1.5)  node{$\scriptstyle{\frac{1-q^1t^4}{\textcolor{red}{1-q^2t^4}}}$};
				\draw (2.5,1.5)  node{$\scriptstyle{\frac{\textcolor{red}{1-q^2t^4}}{1-q^1t^3}}$};
				\draw (3.5,1.5)  node{$\scriptstyle{\frac{\textcolor{red}{1-q^3t^4}}{1-q^0t^2}}$};
				\draw (0.5,2.5)  node{$\scriptstyle{\frac{1-q^0t^3}{\textcolor{red}{1-q^3t^3}}}$};
				\draw (1.5,2.5)  node{$\scriptstyle{\frac{1-q^1t^3}{\textcolor{red}{1-q^2t^3}}}$};
				\draw (2.5,2.5)  node{$\scriptstyle{\frac{\textcolor{red}{1-q^2t^3}}{1-q^1t^2}}$};
				\draw (3.5,2.5)  node{$\scriptstyle{\frac{\textcolor{red}{1-q^3t^3}}{1-q^0t^1}}$};
				\draw (0.5,3.5)  node{$\scriptstyle{\frac{1-q^0t^2}{\textcolor{red}{1-q^2t^2}}}$};
				\draw (1.5,3.5)  node{$\scriptstyle{\frac{\textcolor{red}{1-q^1t^2}}{\textcolor{red}{1-q^1t^2}}}$};
				\draw (2.5,3.5)  node{$\scriptstyle{\frac{\textcolor{red}{1-q^2t^2}}{1-q^0t^1}}$};
				\draw (0.5,4.5) node{$\scriptstyle{\frac{\textcolor{red}{1-q^0t^1}}{\textcolor{red}{1-q^1t^1}}}$};
				\draw (1.5,4.5) node{$\scriptstyle{\frac{\textcolor{red}{1-q^1t^1}}{\textcolor{red}{1-q^0t^1}}}$};
				
			\end{tikzpicture}
		\end{center}

    \caption{Adding missing factors. The product of all entries in the diagram remains equal to the RHS of (\ref{elliptic}).}
    \label{fig:addmissingfactors}
\end{figure}		

		The wonderful thing is that the added red factors in the numerators and denominators are exactly the same. So we aren't actually adding anything new to the product, only rewriting it in a neater way.
		
		Thus, this means that the product of the factors we originally computed from the right-hand side of (\ref{elliptic}) is indeed equal to the left-hand side of the equation. So, this completes our verification of the elliptic version of the theorem for our example.
		
\section{The last upgrade: an introduction to Macdonald polynomials.}
		
		So far, we have introduced the integer version, the polynomial version and the elliptic version of the theorem, with each more general than the one before. A natural question to ask is --- can there be an even more general form?
		
		\textit{Symmetric polynomials} are polynomials that are invariant when interchanging variables \cite[Ch.\ I §2]{MR1354144}. For example, the function $e_2$ below is symmetric.		\[e_2(x_1,x_2, x_3) = x_1x_2+x_1x_3+x_2x_3.\] The function $e_2$ belongs to a special family of symmetric functions called the elementary symmetric functions. In addition, there are many other noteworthy examples including the monomial symmetric functions, the Schur functions, the Hall–Littlewood symmetric functions, etc.

		In 1988, Ian G. Macdonald introduced the symmetric Macdonald polynomials $P_\lambda(x_1, x_2, \ldots, x_n ;q,t)$ indexed by partitions $\lambda$ and two parameters $q, t$ \cite{Macdonald88anew}.

		The Macdonald polynomials unify many symmetric functions (see \cite[Ch.\ I §11]{MR1488699}). For example, when $t=1$, the Macdonald polynomials collapse into the monomial symmetric functions. Further, if we substitute $q=t$, the Macdonald polynomials become the Schur functions. Now, let us consider the principal specialisation of the Macdonald polynomials $P_\lambda(1,t,t^2,\ldots, t^{n-1};q,t)$, which we denote as $P_\lambda(q,t)$ from now on. Interestingly, it turns out that this special case gives us the right hand side of the elliptic version of the theorem (\ref{elliptic}).
		
		\begin{equation*}
				\prod_{b\in \lambda} \dfrac{1-q^{coarm_\lambda(b)}t^{n-coleg_\lambda(b)}}{1-q^{arm_\lambda(b)}t^{leg_\lambda(b)+1}} 
				=
				\prod_{1\leq i < j \leq n} \prod_{r=0}^{\lambda_i-\lambda_j-1} \dfrac{1-q^rt^{j-i+1}}{1-q^rt^{j-i}}   \quad\hbox{comes from}\quad P_\lambda(q,t).
		\end{equation*}
		
		Recall that setting $q=t$ simplifies the elliptic version to the polynomial version of the theorem. Meanwhile, the same substitution $q=t$ reduces $P_\lambda(q,t)$ to the Schur functions $s_\lambda$. This means our polynomial version of the theorem  (\ref{polynomial}) corresponds to the principal specialisation of the Schur functions.
		
		\begin{equation*}
				\prod_{b\in \lambda} \dfrac{1-t^{n+c(b)}}{1-t^{h(b)}} = \prod_{1\leq i <j \leq n} \dfrac{1-t^{\lambda_i-\lambda_j+j-i}}{1-t^{j-i}}  \quad\hbox{comes from }\quad P_\lambda(t,t) = s_\lambda(1, t. \ldots,t^{n-1}).
		\end{equation*}
		
\begin{figure}[h]
    \centering
		\begin{center}
			\begin{tikzpicture}[yscale=1, xscale=1]
				\draw [-stealth](0,0) -- (6,0);
				\draw [-stealth](0,0) -- (0,6);
				\draw (6.5,0)  node{$q$};
				\draw (0,6.5)  node{$t$};
				\filldraw [draw=black, fill=blue, opacity=0.1] (0,0) rectangle (5,5);
				\draw (5,-0.3)  node{$1$};
				\draw (-0.3,5)  node{$1$};
				\filldraw (0.7,4) circle (1pt);
				\draw (1.5,4)  node{$P_\lambda(q,t)$};
				\draw (1.5,3.4)[color=gray, text width=2.8cm, align=center]  node{Elliptic version theorem};
				\draw (0,0) -- (5,5);
				\filldraw (2,2) circle (1pt);
				\draw (3.3,2)  node{$s_\lambda = P_\lambda(t,t)$};
				\draw (3.3,1.3)[color=gray, text width=2.8cm, align=center]  node{Polynomial version theorem};
				\draw[blue, thick, -latex] (4,4) -- (5,5);
				\filldraw (5,5) circle (1pt);
				\draw (5.3,5.3)  node{$\lim_{t\rightarrow 1}P_\lambda(t,t)$};
				\draw (6.3,4.6)[color=gray, text width=2.8cm, align=center]  node{Integer version theorem};
			\end{tikzpicture}
		\end{center}
    \caption{Visualising the three levels of generalisations.}
    \label{fig:square3levels}
\end{figure}
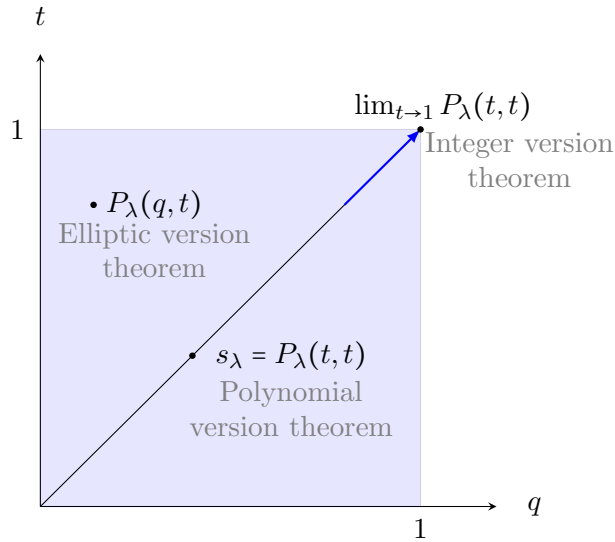

		The values of $q$ and $t$ can vary in countless ways. Let us visualise it using a graph. Setting the horizontal and the vertical axis as $q$-axis and $t$-axis, consider the unit square at the origin. Each point in this square represents a family of symmetric polynomials in the elliptic version $P_\lambda(q,t)$. Then, the line $t=q$ corresponds to the polynomial version of the theorem (\ref{polynomial}) and the limit as we approach the point $(1,1)$ along the line $t=q$ represents the integer version (\ref{integer}) as indicated in Figure \ref{fig:square3levels}.

\begin{figure}[h]
    \centering
    		\begin{center}
			\begin{tikzpicture}[yscale=1, xscale=1]
				\draw [-stealth](0,0) -- (6,0);
				\draw [-stealth](0,0) -- (0,6);
				\draw (6.5,0)  node{$q$};
				\draw (0,6.5)  node{$t$};
				\filldraw [draw=black, fill=blue, opacity=0.1] (0,0) rectangle (5,5);
				\draw (5,-0.3)  node{$1$};
				\draw (-0.3,5)  node{$1$};
				\filldraw (0.7,4) circle (1pt);
				\draw (1.5,4)  node{$P_\lambda(q,t)$};
				\draw[color=gray, opacity = 0.8, text width = 2.2cm, align = center] (1.5,3.3)  node{Macdonald polynomials};
				\draw (0,0) -- (5,5);
				\filldraw (2,2) circle (1pt);
				\draw (3.3,2)  node{$s_\lambda = P_\lambda(t,t)$};
				\draw[color=gray, opacity = 0.8] (3.3,1.6)  node{Schur};
				\draw[color=gray, opacity = 0.8] (3.3,1.3)  node{functions};
				
				\draw[color=blue] (0,5) -- (5,5);
				\filldraw (2,5) circle (1pt);
				\draw[color=black] (2,5.3)  node{$m_\lambda = P_\lambda(q,1)$};
				\draw[color=gray, opacity = 0.8] (2,5.7)  node{symmetric functions};
				\draw[color=gray, opacity = 0.8] (2,6.1)  node{Monomial};
				 
				\draw[color=blue] (5,0) -- (5,5);
				\filldraw (5,2.5) circle (1pt);
				\draw[color=black] (6.2,2.6)  node{$e_{\lambda'} = P_\lambda(1,t)$};
				\draw[color=gray, opacity = 0.8] (6.2,2.2)  node{Elementary};
				\draw[color=gray, opacity = 0.8] (6.2,1.8)  node{symmetric};
				\draw[color=gray, opacity = 0.8] (6.2,1.4)  node{functions};
				 
				\draw[color=blue] (0,0) -- (0,5);
				\filldraw (0,2) circle (1pt);
				\draw[color=black] (-0.8,2)  node{$P_\lambda(0,t)$};
				\draw[color=gray, opacity = 0.8] (-1.4,1.6)  node{Hall-Littlewood};
				\draw[color=gray, opacity = 0.8] (-1.4,1.2)  node{polynomials};
				 
				\draw[color=blue] (0,0) -- (5,0);
				\filldraw (2,0) circle (1pt);
				\draw[color=black] (2,-0.3)  node{$P_\lambda(q,0)$};
				\draw[color=gray, opacity = 0.8] (2,-0.7) node{q-Whittaker polynomials};
				 
				\draw[blue, thick, latex-] (5,5) arc (350:285:1cm);
				\draw[blue, thick, latex-] (5,5) arc (355:300:1.2cm);
				\draw[blue, thick, latex-] (5,5) arc (358:330:2.5cm);
				\filldraw (5,5) circle (1pt);
				\draw[color=black] (5.5,5.5)  node{$P_\lambda^{(\alpha)}=\lim\limits_{t\rightarrow 1} P_\lambda(t^\alpha,t)$};
				\draw[color=gray, opacity = 0.8] (5.6,5) node{Jack};
				\draw[color=gray, opacity = 0.8] (6,4.6) node{symmetric};
				\draw[color=gray, opacity = 0.8] (6,4.2) node{functions};
			\end{tikzpicture}
		\end{center}
    \caption{The principal specialisation square for the Macdonald polynomials.}
    \label{fig:square}
\end{figure}
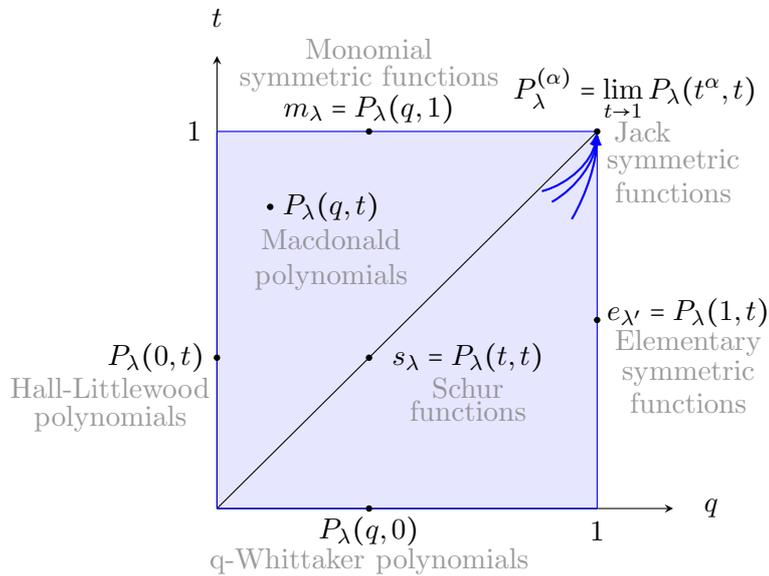

		We mentioned earlier that there are also some other types of symmetric functions --- they can also be added to the picture, as shown in Figure \ref{fig:square}!
        \begin{itemize}
    	    \item When $t=1$, $P_\lambda(q,1)=m_\lambda$ is the monomial symmetric function.
    	    \item When $q=1$, $P_\lambda(1,t)=e_{\lambda'} $ is the elementary symmetric function.
    	    \item When $q=0$, $P_\lambda(0,t)$ is the Hall-Littlewood polynomial.
    	    \item When $t=0$, $P_\lambda(q,0)$ is the $q$-Whittaker polynomial.
    	    \item When $q=t^\alpha$ and $t \rightarrow 1$, $P_\lambda^{(\alpha)}$ is the Jack polynomial. 
        \end{itemize}		

		Isn't that magnificent? It all started with something as simple as counting boxes.
		
		\nocite{*}
\section{The becoming of this article}
    	`I have two tasks on my desk that I need help with $\ldots$' That was how Professor Arun Ram lured us into this fascinating world of the Macdonald polynomials and we wouldn't have it any other way. The three of us are undergraduate students from the University of Melbourne. We met each other in Arun's Group Theory and Linear Algebra class and have been exploring interesting enrichment topics with him ever since. From our first time seeing these elegant proofs being presented by Arun during Zoom lectures, to discussing ideas among ourselves and getting so excited by our combinatorial interpretation of the formulas, to then presenting at our student maths society seminar to test the waters of our content --- the becoming of this article has been such a rewarding, fun and intellectually stimulating journey. We hope you enjoyed reading it as much as we enjoyed writing it!
    	
\section{Acknowledgements}
        We would like to express our deepest thanks to Professor Arun Ram for introducing us to the interesting topics of Young tableaux and Macdonald polynomials,  and for his patience and generosity in guiding three unsuspecting undergraduate students through their first article.
        
        All three of us would also like to thank each other for the useful discussions and the team spirit of helping each other throughout the process. The enthusiasm for mathematics from each one of us inspires the whole group to keep doing maths with enjoyment. 
		\printbibliography

\end{document}